\documentclass{article}
\usepackage{graphicx} % Required for inserting images
\usepackage{hyperref} 
\title{There are more non-Cantorian sets than are Cantorian}
\author{Zuhair Al-Johar}
\date{March 2025}

\begin{document}

\maketitle

\section{Introduction}

When working in $\sf NF$, [1] there is a sense that there are more non-Cantorian sets than Cantorian sets. But it is not that immediate result as one expects, since they are externally equinumerous, and the qualification ``Cantorian" is not stratified and so not easy to spell internally. This account stipulates a fairly simple criterion to phrase such problems and proves that per that criterion there are more non-Cantorian sets.

\section{A criterion of comparison of predicate sizes}
We say that a predicate $P$ has strictly fewer objects that fulfill it than those that fulfill a predicate $Q$, denoted as $P \prec Q$, if and only if there is a set $X$ of which every object fulfilling $P$ is an element, and a set $Y$ that has all its elements fulfilling $Q$, and $|X| < |Y|$. 

The operator $|X|$ is defined, as usual in $\sf NF$, after Frege, as the set of all sets that are bijective to $X$, [2] and the operator ``$<$" is defined, as usual, by Cantor, as the existence of an injective function from left to right and absence of it in the opposite direction. [3]

\section {The claim and its proof} 

\textbf{Claim:} Cantorian $\prec$ non-Cantorian \smallskip

Where a set $x$ is said to be Cantorian, if and only if, $|x| =|\iota``x|$, where $\iota``x=\{\{y\} : y \in x\}$. [4] \bigskip

\textbf{Proof: } Let $X$ be the set of all sets that are surjective images of $\iota``V$ (i.e. there exist surjections from $\iota``V$ onto each of them), where $V$ is the set of all sets. Formally, $X=\{S : |S| \leq^* |\iota``V| \}$.

Let $Y = \{ V \times x: x \in V \}$, where the operator ``$\times$" stands for the Cartesian product that uses the Quine-Rosser pairs. [5] \smallskip

Every Cantorian set $s$ would be an element of  $X$, since $s$ is injective to $\iota``V$ and so $s$ would be a surjective image of $\iota``V$.\smallskip

All elements of $Y$ are non-Cantorian since they are $V$-sized, this is directly shown from Schröder–Bernstein theorem ``$\sf SB$", [6] which $\sf NF$ proves, [4] and the maps $x \mapsto \langle 0,x\rangle; x \mapsto x$. Also we note that $Y$ itself is $V$-sized by the maps $x \mapsto x \times V; x \mapsto x$.   \bigskip

The proof that $|X| < |Y|$: 

Bowler in [7]  has shown that the set of all surjective images of $\iota``V$, is itself a surjective image of $\iota``V$, so we have $X$ is a surjective image of $\iota``V$. Bowler's proof involves the construction of a map: $\phi: \iota``V \to V$, defined as: $$\phi: \{C\} \mapsto \{\{s: \langle x,s \rangle \in C\}: x \in V\}$$ So, each $\phi(\{C\})$ is a surjective image of $\iota``V$, and each $S$ whereby there exists a surjection $F:\iota``V \to S$, we take $C=\{\langle x,s \rangle : s \in F(\{x\})\}$, and we'll have $S=\phi(\{C\})$. Consequently, $X$ is the image of $\phi$, which is a surjection from $\iota``V$, that is, $X=\{S \mid \exists x \in \iota``V: \langle x,S \rangle \in \phi \}$.\bigskip

Now, $X$ is injective to $Y$, since every set is injective to $Y$ because $Y$ is $V$-sized.
However, $Y$ cannot be injective to $X$, since otherwise this would mean that there exists a bijection between $X,Y$ by 
$\sf SB$, [6] but this would mean that $V $ is a surjective image of $\iota``V$. But that cannot be since it's known that if $ B$ is a surjective image of $A$, then $\mathcal P(B)$ injects to $\mathcal P(A)$, accordingly we will have $\mathcal P(V)$ injecting into $\mathcal P(\iota``V)$, and so $V$ injects into $\iota``V$ since $|\mathcal P(\iota``V)|=|\iota``V|$ by maps $x \mapsto\{\bigcup x\}; \{x\} \mapsto  \iota``x$; but by then $V$ would be Cantorian, which is not! \bigskip

\section{References} 
\begin{enumerate}
    \item Quine WV. New Foundations for Mathematical Logic. Am Math Mon. 1937;44:70–80.
    \item Frege G. The Foundations of Arithmetic (Die Grundlagen der Arithmetik, with facing English text by J. L. Austin). Evanston: Northwestern University Press; 1976.
    \item Cantor G. Beiträge zur Begründung der transfiniten Mengenlehre (Contributions to the Founding of the Theory of Transfinite Numbers, 1895–1897).
    \item Holmes MR. Elementary set theory with a universal set [PDF]. Cahiers du Centre de Logique. Vol. 10. Louvain-la-Neuve: Université Catholique de Louvain, Département de Philosophie; 1998. ISBN: 2-87209-488-1.
    \item Rosser JB. Logic for Mathematicians. New York: McGraw-Hill; 1953.
    \item Bernstein F. Untersuchungen aus der Mengenlehre. Math Ann. 1905;61:117–155.
    \item Forster T. Scrapbook on Set Theory with a Universal Set. Preprint. 2025 Jan 1. p. 310.
    Available from:\url{https://www.dpmms.cam.ac.uk/\~tef10/NFnotesredux.pdf}. Accessed 2025 Mar 7.
\end{enumerate}

\end{document}